\documentstyle[12pt,fullpage]{article}

\begin{document}

\renewcommand{\thesubsection}{\arabic{subsection}} 
\newtheorem{theorem}{Theorem} 
\newtheorem{prop}[theorem]{Proposition}
\newtheorem{lemma}[theorem]{Lemma}
\newtheorem{claim}[theorem]{Claim}
\newtheorem{corollary}[theorem]{Corollary}
\newtheorem{remark}[theorem]{Remark}
\newtheorem{conjecture}[theorem]{Conjecture}
\newtheorem{example}[theorem]{Example}
\newtheorem{obsr}[theorem]{Observation}
\newtheorem{fact}[theorem]{Fact}
\newtheorem{defn}[theorem]{Definition}
\newtheorem{exercise}[theorem]{Exercise}
\newtheorem{algo}{Algorithm}[section]

\newlength{\originalbase}
\setlength{\originalbase}{\baselineskip}
\newcommand{\spacing}[1]{\setlength{\baselineskip}{#1\originalbase}}

\newenvironment{proof}{\noindent{\bf Proof}:\,\,} {\hfill $\Box$\medskip}

\def\un{{\mathtt Df}}
\def\df{{\mathtt df}}
\def\odd{{\mathtt od}}
\def\Odd{{\mathtt Odd}}
\def\ev{{\mathtt ev}}
\def\Ev{{\mathtt Ev}}
\def\bsl{{\,\!-\!\,}}
\def\T{{\sf ^T}}
\def\R{\hbox{\sf I\kern-.12em R}}
\def\sR{\hbox{\scriptsize\sf I\kern-.12em R}}
\def\N{\hbox{\sf I\kern-.12em N}}
\def\Q{\hbox{\sf C\kern-.47em Q}}
\def\C{\hbox{\sf C\kern-.47em C}}
\def\sC{\hbox{\scriptsize\sf C\kern-.47em C}}
\def\mid{\mbox{$\;:\;$}}
\def\proofend{\hfill$\Box$

\medskip}

\begin{center}
{\Large\sf
Short Proof of the Gallai-Edmonds Structure Theorem
}\\[8mm]
{\sc Andre\u{\i} Kotlov}\\[2mm]
{\small CWI, Kruislaan 413, 1098 SJ Amsterdam, the Netherlands}\\
{\tt andrei$@$cwi.nl}\\[6mm]
\end{center}

\begin{abstract}
{\small
\baselineskip=0.158in
We derive the
Gallai-Edmonds Structure Theorem from Hall's Theorem.
}
\end{abstract}

\bigskip

\noindent{\large\bf Introduction and Definitions}

\bigskip

\noindent
In 1971, Anderson~\cite{a} gave a derivation of Tutte's Theorem from
Hall's Theorem. In this note, we give a concise proof of the
Gallai-Edmonds Structure Theorem; our proof-method 
seems new but similar to that of Anderson.  For completeness, we
proceed with the basic definitions and relevant formulations;  the
reader familiar with the Gallai-Edmonds Structure Theorem (see, for
example, \cite{lp} pp.93--95) is
asked to skip to the next section.

A set $M$ of edges in a graph $G$ is {\bf a matching} if no two edges
of $M$ share a vertex.  A vertex of $G$  incident with no edge in 
$M$ is {\bf $M$-exposed}.  A matching  is {\bf perfect} if it has no exposed
vertex, and {\bf near-perfect} if it has exactly one exposed vertex.
A graph $G$ is {\bf factor-critical} if for every vertex $v\in G$ the
graph $G\bsl v$ has a perfect matching.   

In what follows, $S$ and $T$ are subsets of the vertex set of
$G$; we set $s:=|S|$ and $t:=|T|$.  
{\bf The [odd, even] $S$-components} are the connected components of
the graph $G\bsl S$ [of odd, even cardinality, respectively].  The number of
odd $S$-components is denoted by $\odd(S)$, the total number of
vertices in those $S$-components which have no perfect matching
is denoted by $\boldmath\un(S)$, while the quantity
$\df(S):=\odd(S)-s$ is {\bf the  
deficiency} of $S$.  It is important to notice that (1) every matching in
$G$ leaves at least $\df(S)$ nodes exposed, and (2) if a matching $M$
leaves exactly $\df(S)$ nodes exposed then $M$ contains  a perfect
matching of each even $S$-component, a near-perfect matching of each
odd $S$-component, and matches all nodes of $S$ to nodes in distinct 
$S$-components.  If (2) is the case for some $S$ and $M$, then $S$ is
called {\bf Tutte-Berge}, while $M$ is maximum by (1).  {\bf The Tutte-Berge
Formula} asserts that a Tutte-Berge set $S$ always exists.  

%

If $s+\odd(S)\ge 2$, we denote by $\langle G,S\rangle$ the bipartite
{\em minor} of $G$ obtained from $G$ by deleting the vertices of the
even $S$-components, contracting each odd 
$S$-component to a single node, and deleting the edges spanned by
$S$.  A monochromatic set $S$ of vertices in a bipartite graph $H$ {\bf
satisfies Hall's condition [{\rm resp.,} with surplus $k\in\N$]} if 
for every non-empty subset $T$ of $S$, the size of the neighborhood of 
$T$ in $H$ is at least $t$ [resp., $t+k$].  In these terms, {\bf Hall's
Marriage Theorem} asserts that {\em $S$ is covered in $H$ by some matching
if and only if $S$ satisfies Hall's condition.} 
Finally, $S$ is {\bf Gallai-Edmonds [with respect to $G$]} if:

(a) {the even $S$-components, if any, have a perfect matching;}

(b) {the odd $S$-components, if any, are factor critical;}

(c) {if $S$ is non-empty, then $S$ satisfies Hall's condition with
surplus one in $\langle G,S\rangle$.} 

\bigskip

\bigskip

\noindent{\large\bf The Gallai-Edmonds Structure Theorem}

\medskip

\begin{itemize}
\item[(i)] {\em For every graph $G$, there exists a Gallai-Edmonds set,
$S$.}

\item[(ii)] {\em $S$ is Tutte-Berge. Consequently, every maximum
matching of $G$ contains a near-perfect matching of each odd
$S$-component, a perfect matching of each  
even $S$-component, and matches all nodes of $S$ to nodes in distinct
$S$-components.}

\item[(iii)] 
{\em The underlying vertex set of the odd $S$-components is the set 
$D(G)$ of the vertices left exposed by at least one maximum matching
of $G$, while $S$ is the neighborhood of $D(G)$.  In particular, $G$
has a unique Gallai-Edmonds set.}
\end{itemize}
\medskip

\noindent{\bf Proof:\,}
Among the subsets of the vertex set of $G$ with {\em maximum} deficiency,
let $S$ have {\em minimum} $\un(S)$.  We show, by induction on $|G|$, 
that $S$ is Gallai-Edmonds.  

Suppose that $C$ is an even $S$-component 
with no
perfect matching.  Fix $v\in C$ and set 
$S':=S\cup\{v\}$.  Then $\df(S')\ge\df(S)$ while $\un(S')<\un(S)$.
This contradiction shows that $S$ satisfies (a).

Suppose that $v$ is a vertex in an odd $S$-component $C$ such
that $C\bsl v$ has no perfect matching.  By induction, the graph
$H:=C\bsl v$ has a Gallai-Edmonds set, $T$.  In particular,
$\df_H(T)\ge 2$. Set $S':=S\cup T\cup\{v\}$. 
Then $\df(S')=\odd(S')-|S'|=[\odd(S)+\odd_H(
T)-1]-[s+t+1]=\df(S)+\df_H(T)-2\ge\df(S)$ while $\un(S')<\un(S)$.  
Thus, $S$ satisfies (b). 

Suppose that $T$ is a smallest non-empty subset of $S$ violating
Hall's condition {\em with surplus one} in  $\langle G,S\rangle$.  Set
$S':=S\bsl T$.  If $T$ consists of a single vertex 
with no neighbors in $\langle G,S\rangle$ then $\df(S')>\df(S)$ which 
is a contradiction.  Else, $T$ satisfies Hall's condition in $\langle
G,S\rangle$. By (b) and Hall's Theorem, $T$ is contained in an even
$S'$-component 
with a perfect matching (this even component is spanned by $T$ and the $t$
odd $S$-components ``neighboring'' with $T$ in $G$).  Consequently,
$\df(S')=\df(S)$ and $\un(S')<\un(S)$.  Thus, $S$ satisfies (c) whence 
(i). 

Let now $S$ be Gallai-Edmonds, and let $v$ be a fixed vertex in an odd
$S$-component $C$. By (c) and Hall's Theorem, $S$ can be matched to
nodes in $s$ distinct odd 
$S$-components different from $C$. By (a) and (b), this matching can be
extended by (near-)perfect matchings of the
$S$-components to obtain a matching $M$ avoiding $v$.  Since $M$
leaves $\df(S)$ nodes exposed, it is maximum, and (ii) and (iii)
follow immediately.  
\proofend

\bigskip

\bigskip

\noindent{\large\bf Concluding Remarks}

\medskip

\noindent\begin{enumerate}
\item{The above proof can be formulated as an exercise, as follows: {\em (A) 
Among the sets with maximum deficiency, consider a set $S$ which
minimizes the total number of vertices in the $S$-components with no
perfect matching; (B) prove, by contradiction, that $S$ is
Gallai-Edmonds; (C) deduce the Gallai-Edmonds Structure Theorem.}}

\item{Among the sets $T$ with minimum number of even
$T$-components with no perfect matching, we could choose $S$ to maximize
the difference between $\df(S)$ and the total number of vertices in the 
odd $S$-components. The above proof would then be repeated almost verbatim.}

\item{The proof is further shortened if the Tutte-Berge Formula is
taken for granted.  In fact, $S$ can be chosen among the Tutte-Berge
sets to minimize the total number of vertices in the  odd
$S$-components. Then (a) and (ii) become superfluent for $S$, while
(b) and (c) can be shown as above but simpler. }

\item{By incorporating the proof of Hall's Theorem into the above
inductive argument, one can derive the Gallai-Edmonds Structure
Theorem ``from scratch.''}
\end{enumerate}

\bigskip

\noindent{\large \bf Acknowledgements}

\bigskip

\noindent
I would like to express my gratitude to Noga Alon, 
 Romeo Rizzi, and especially Bert Gerards and Lex
Schrijver, for their indispensable help.  
%

%

\end{document}